\documentclass{amsart}
\usepackage{amssymb,amsthm,amstext,amsmath,bm}
\usepackage[inline]{enumitem}

\usepackage{hyperref} 

\usepackage{xcolor} 

\newcommand{\N}{\ensuremath{\mathbb{N}}}
\newcommand{\Q}{\ensuremath{\mathbb{Q}}}
\newcommand{\E}{\mathfrak{E}}
\newcommand{\Ec}{\mathfrak{E}_c}
\newcommand{\cl}[2][X]{\mathrm{cl}_{#1}\!\left(#2\right)}
\newcommand{\ext}[2][C]{\mathrm{ext}_{#1}(#2)}
\newcommand{\vietoris}[1]{{\boldsymbol{\langle\!\!\langle}}#1 {\boldsymbol{\rangle\!\!\rangle}}}

\newtheoremstyle{theorem}
     {11pt}
     {11pt}
     {}
     {}
     {\bfseries}
     {}
     {.5em}
     {\noindent\thmnumber{#2}. \thmname{#1}{\rm\thmnote{#3}}}

\theoremstyle{theorem}

\newtheorem{defi}{Definition}[section]
\newtheorem{theorem}[defi]{Theorem}
\newtheorem{propo}[defi]{Proposition}
\newtheorem{lemma}[defi]{Lemma}
\newtheorem{cor}[defi]{Corollary}
\newtheorem{ques}[defi]{Question}

\newtheorem{remark}[defi]{Remark}
\newtheorem{ex}[defi]{Example}

\title[A characterization of $\Q\times\Ec$]{A characterization of the product of the rational numbers and complete Erd\H{o}s space}
\author[R. Hernández-Gutiérrez]{Rodrigo Hernández-Gutiérrez}
\author[A. Zaragoza]{Alfredo Zaragoza}
\address[R. Hernández-Gutiérrez]{Departamento de Matem\'aticas, Universidad Aut\'onoma Metropolitana campus Iztapalapa, Av. San Rafael Atlixco 186, Col. Vicentina, Iztapalapa, 09340, Mexico city, Mexico}
\email[R. Hernández-Gutiérrez]{rod@xanum.uam.mx}
\address[A. Zaragoza]{Departamento de Matemáticas, Facultad de Ciencias, Universidad Nacional Autónoma de México, Circuito Exterior s/n, Ciudad Universitaria, Coyoacán, 04510, Mexico city, Mexico
}
\email[A. Zaragoza]{soad151192@icloud.com}
\thanks{This work is part of the doctoral work of the second-named author at UNAM, Mexico city, under the direction of the first-named author. This research was supported by a CONACyT doctoral scholarship with number 696239.}

\date{\today}


\begin{document}

\begin{abstract}
Erdős space $\mathfrak{E}$ and complete Erdős space $\mathfrak{E}_c$ have been previously shown to have topological characterizations. In this paper, we provide a topological characterization of the topological space $\mathbb{Q}\times\mathfrak{E}_c$, where $\mathbb{Q}$ is the space of rational numbers. As a corollary, we show that the Vietoris hyperspace of finite sets $\mathcal{F}(\mathfrak{E}_c)$ is homeomorphic to $\mathbb{Q}\times\mathfrak{E}_c$. We also characterize the factors of $\mathbb{Q}\times\mathfrak{E}_c$. An interesting open question that is left open is whether $\sigma\mathfrak{E}_c^\omega$, the $\sigma$-product of countably many copies of $\mathfrak{E}_c$, is homeomorphic to $\mathbb{Q}\times\mathfrak{E}_c$.
\end{abstract}

\maketitle

\section{Introduction}

All spaces will be assumed to be separable and metrizable. We denote the set of positive integers by $\N$, the set of natural numbers by $\omega=\N\cup\{0\}$ and the space of rational numbers by $\Q$. Erdős space is defined to be the space
$$ \E = \{(x_n)_{n\in \omega} \in  \ell^2 \colon \forall i\in\omega,\, x_i \in  \Q\},$$
and complete Erdős space is the space
$$ \Ec = \{(x_n)_{n\in \omega} \in  \ell^2 \colon \forall i\in\omega,\, x_i \in\{0\}\cup\{1/n:n\in\N\}\},$$
where $\ell^2$ is the Hilbert space of square-summable sequences of real numbers.
These two spaces were introduced by Erdős in 1940 in \cite{er} as examples of totally disconnected and non-zero-dimensional spaces.

It was soon noticed that some interesting Polish spaces are homeomorphic to $\Ec$, see \cite{k-o-t}. Due to the interest in these two spaces, Jan Dijkstra and Jan van Mill obtained topological characterizations of $\Ec$ and $\E$ (see \cite{DvM} and \cite{ME}, respectively), and applied them to show that  some other noteworthy spaces are homeomorphic to either one of these two. Notice that $\Ec$ is Polish but $\E$ is an absolute $F_{\sigma\delta}$, so $\Ec$ and $\E$ are not homeomorphic. We also mention that $\Ec^\omega$ is not homeomorphic to $\Ec$, as it was proved in \cite{d-vm-s}. A characterization of $\Ec^\omega$ was given in \cite{dijkstra-stable}.

The objective of this paper is to continue this line of research by providing a topological characterization of $\Q\times\Ec$; this is Theorem \ref{thm-main} below. Since $\Q\times\Ec$ is not Polish, it is not homeomorphic to $\Ec$ or $\Ec^\omega$. As it is easy to see, $\Q\times\Ec$ is both an absolute $G_{\delta\sigma}$ and an absolute $F_{\sigma\delta}$ (see Remark \ref{complexity-QEc} below). Since it is known that $\E$ is not $G_{\delta\sigma}$ (see Remark 5.5 in \cite{ME}), we obtain that $\Q\times\Ec$ is not homeomorphic to $\E$. Thus, this space is different from the ones studied before.

In fact, we give two characterizations of $\Q\times\Ec$\,: one extrinsic and the other intrinsic. The choice of these two terms follows the idea of \cite{ME}. By extrinsic we mean that $\Q\times\Ec$ is homeomorphic to a subset of the graph of a USC (upper semi-continuous) function defined on the Cantor set that has certain characteristics. By intrinsic we mean a characterization given by topological properties of $\Q\times\Ec$ itself. Our extrinsic characterization is defined in terms of a class $\sigma\mathcal{L}$ of USC functions and our intrinsic characterization is given by a class $\sigma\mathcal{E}$ of spaces; both of these are defined in section \ref{section:classes}. 

The statement of the characterization, Theorem \ref{thm-main}, is given in section 3 but the hard part of the proof is done in section \ref{section:main-theorem}. We also give a concrete application of our characterizations: in section \ref{section:hyperspace} the Vietoris hyperspace of finite non-empty subsets of $\Ec$ is shown to be homeomorphic to $\Q\times\Ec$ (Corollary \ref{hyperspace}). This result is connected to previous work of the second-named author who proved that the Vietoris symmetric products of $\Ec$ are homeomorphic to $\Ec$ (see \cite{zaragoza-1}) and that the Vietoris hyperspace of non-empty finite sets of $\E$ is homeomorphic to $\E$ (see \cite{zaragoza-1} and \cite{zaragoza-2}). In section \ref{section:sigma-product} we consider the $\sigma$-product of $\omega$ copies of $\Ec$. At first, it seemed that this space would also be homeomorphic to $\Q\times\Ec$. However, we were not able to prove or disprove this, so we leave this as an open problem. In section \ref{section:factors} we give a characterization of factors of $\Q\times\Ec$. Finally, in section \ref{section:embeddings} we consider dense embeddings of $\Q\times\Ec$.

\section{Preliminaries}

Following the example of E. K. van Douwen, we call a space \emph{crowded} if it has no isolated points. The definitions and equivalences that we will use here can be found in \cite{ME}. The notation $X\approx Y$ means that $X$ and $Y$ are homeomorphic topological spaces.

A \emph{$C$-set} in a topological space is an intersection of clopen sets. A topological space is \emph{almost zero-dimensional} if it has a neighborhood basis consisting of $C$-sets. Given a topological space $\langle X,\mathcal{T}\rangle$ and $A\subset X$ we write $\mathcal{T}\restriction A=\{U\cap A\colon U\in\mathcal{T}\}$.

\begin{defi}
 Let $\langle X,\mathcal{T}\rangle$ be a topological space and let $\langle Z,\mathcal{W}\rangle$ be a zero-dimensional space such that $X\subset Z$. We will say that $\langle Z,\mathcal{W}\rangle$ \emph{witnesses the almost zero-dimensionality} of  $\langle X,\mathcal{T}\rangle$ if $\mathcal{W}\restriction X\subset\mathcal{T}$ and there is a neighborhood basis of $\langle X,\mathcal{T}\rangle$ that consists of sets that are closed in $\mathcal{W}$.
\end{defi}

It easily follows that a topological space $\langle X,\mathcal{T}\rangle$ is almost zero-dimensional if and only if there is a zero-dimensional topology $\mathcal{W}$ in $X$ that witnesses the almost zero-dimensionality of $\langle X,\mathcal{T}\rangle$, see \cite[Remark 2.4]{ME}.

Let $X$ be a space and let $\mathcal{A}$ be a collection of subsets of $X$. The space $X$ is called \emph{$\mathcal{A}$-cohesive} if every point of the space has a neighborhood that does not contain non-empty clopen subsets of any element of $\mathcal{A}$. If $\mathcal{A}=\{X\}$, we simply say that $X$ is cohesive.

Let $\varphi\colon X\to[0,\infty)$. We say that $\varphi$ is USC (upper semi-continuous) if for every $t\in(0,\infty)$ the set $f\sp\leftarrow[(-\infty,t)]$ is open. Let
$$
M(\varphi)=\sup\left(\{\lvert\varphi(x)\rvert\colon x\in X\}\cup\{0\}\right),
$$
where the supremmum is taken in $[0,\infty]$. We define
 $$
 \begin{array}{ccl}
  G_0^\varphi &  = & \{\langle x,\varphi(x)\rangle\colon x\in X,\, \varphi(x)>0\},\ \textrm{and}\\[0.5em]
  L_0^\varphi &  = & \{\langle x,t\rangle\colon x\in X,\, 0\leq t\leq\varphi(x)\}.
 \end{array}
 $$
We say that $\varphi$ is a \emph{Lelek function} if $X$ is zero-dimensional, $\varphi$ is USC, $\{x\in X\colon \varphi(x)>0\}$ is dense in $X$ and $G_0^\varphi$ is dense in $L_0^\varphi$. The existence of Lelek functions with domain equal to the Cantor set $2^\omega$ follows from Lelek's original construction \cite{lelek} of what is now called the Lelek fan.

We will need to extend USC functions. Assume that $X$ is a space, $Y\subset X$ and $\varphi\colon Y\to[0,\infty)$ is a USC function. Then there is a cannonical extension $\ext[X]{\varphi}\colon X\to[0,\infty)$; we will not need its definition (which can be found in \cite[p. 12]{ME}) but only the following property.

\begin{lemma}\cite[Lemma 4.8]{ME}\label{lema:ext-USC}
 Let $X$ be a zero-dimensional space, let $Y$ be a dense subset of $X$, let $\psi\colon Y\to[0,\infty)$ be a USC function and let $\varphi=\ext[X]{\psi}$. Then $\varphi$ is USC, $\psi\subset\varphi$ and the graph of $\psi$ is dense in the graph of $\varphi$.
\end{lemma}

As mentioned in the introduction, $\Ec$ is a cohesive almost zero-dimensional space. An extrinsic characterization of $\Ec$ is given by Lelek functions as follows: if \mbox{$\varphi\colon 2^\omega\to[0,\infty)$} is a Lelek function, then $G_0^\varphi$ is homeomorphic to $\Ec$, see \cite{k-o-t}. An intrinsic characterization of $\Ec$ was given in \cite{DvM}. We make the following remark about the descriptive complexity of $\Q\times\Ec$.

\begin{remark}\label{complexity-QEc}
 $\Q\times\Ec$ is an absolute $G_{\sigma\delta}$ and an absolute $F_{\delta\sigma}$.
\end{remark}
\begin{proof}
 To see that $\Q\times\Ec$ is an absolute $G_{\delta\sigma}$, it is sufficient to notice that $\Q\times\Ec$ is a countable union of Polish spaces.
 
 Next, assume that $\Q\times\Ec\subset X$ where $X$ is any separable metrizable space. For each $q\in\Q$, let $F_q=\{q\}\times\Ec$. Then $G=X\setminus\bigcup\{\overline{F_q}\colon q\in\Q\}$ is a $G_\delta$ in $X$. 
 
 Fix $q\in\Q$. Since $F_q$ is Polish we know that $\overline{F_q}\setminus F_q$ is a countable union of sets that are closed in $\overline{F_q}$, and thus, in $X$. But closed sets in separable metrizable spaces are $G_\delta$. Thus, $\overline{F_q}\setminus F_q$ is $G_{\delta\sigma}$ in $X$.
 
 Since $X\setminus (\Q\times\Ec)=G\cup\left(\bigcup\left\{\overline{F_q}\setminus F_q:q\in\Q\right\}\right)
 $ we obtain that the complement of $\Q\times\Ec$ is a $G_{\delta\sigma}$ so $\Q\times\Ec$ itself is $F_{\sigma\delta}$ in $X$.
\end{proof}

\section{Classes \texorpdfstring{$\sigma\mathcal{L}$}{L} and \texorpdfstring{$\sigma\mathcal{E}$}{sE}}\label{section:classes}

In this section we define the two classes of spaces $\sigma\mathcal{L}$ and $\sigma\mathcal{E}$ that we will use to characterize $\Q\times\Ec$. These definitions are made in the spirit of the class $\mathsf{CAP}(X)$ from \cite{vm-cap} and classes $\mathsf{SLC}$ and $\mathsf{E}$ from \cite{ME}.

\begin{defi}\label{defi:class-L}
We define $\sigma\mathcal{L}$ to be the class of all triples $\langle C, X,\varphi\rangle$ such that $C$ is a compact, zero-dimensional, crowded metrizable space (thus, a Cantor set), $\varphi\colon C\to[0,1)$ is an USC function and $X=\bigcup\{X_n\colon n\in\omega\}$ is a dense subset of $C$ such that for each $n\in\omega$ the following hold
\begin{enumerate}[label=(\alph*)]
    \item $X_n$ is a closed, crowded subset of $C$,
    \item $X_n\subset X_{n+1}$,
    \item $\varphi\restriction{X_n}$ is a Lelek function, and
    \item $G_0^{\varphi\restriction X_n}$ is nowhere dense in $G_0^{\varphi\restriction X_{n+1}}$.
\end{enumerate}
We will say that a space $E$ is generated by  $\langle C, X, \varphi\rangle$ if $E$ is homeomorphic to $G_0^{\varphi \restriction{X}}$.
\end{defi}

As mentioned in the previous section, by the extrinsic characterization of $\Ec$ from \cite{k-o-t}, in Definition \ref{defi:class-L} we will have that $G_0^{\varphi\restriction X_n}$ is homeomorphic to $\Ec$ for each $n\in\omega$. So indeed, $E$ is a countable increasing union of nowhere dense subsets, each homeomorphic to complete Erdős space.

\begin{defi}\label{defi:class-sE}
 We define $\sigma\mathcal{E}$ to be the class of all separable metrizable spaces $E$ such that there exists a topology $\mathcal{W}$ on $E$ that is witness to the almost zero-dimensionality of $E$, a collection $\{E_n\colon n\in\omega\}$ of subsets of $E$ and a basis $\beta$ of neighborhoods of $E$ such that
 \begin{enumerate}[label=(\alph*)]
    \item $E=\bigcup\{E_n\colon n\in\omega\}$,
    \item for each $n\in\omega$,  $E_n$ is a crowded nowhere dense subset of $E_{n+1}$,
    \item for each $n\in\omega$,  $E_n$ is closed in $\mathcal{W}$,
    \item $E$ is $\{E_n\colon n\in\omega\}$-cohesive, and
    \item for each $V\in\beta$, $V\cap E_n$ is compact in $\mathcal{W}\restriction{E_n}$ for each $n\in \omega$.
 \end{enumerate}
\end{defi}

By the intrinsic characterization of $\Ec$ from \cite{DvM}, we have that, in Definition \ref{defi:class-sE}, $E_n$ is homemorphic to $\Ec$ for every $n\in\omega$. So again $E$ is a countable increasing union of nowhere dense subsets, each homeomorphic to complete Erdős space. 

We first prove that the space that we want to characterize is an element of $\sigma\mathcal{E}$ and then, that spaces from $\sigma\mathcal{E}$ can be generated by triples from $\sigma\mathcal{L}$.

\begin{lemma}\label{lema:QtimesEc-sigmaE}
 $\Q\times\Ec\in\sigma\mathcal{E}$.
\end{lemma}
\begin{proof}
By (2) in \cite[Theorem 3.1]{DvM}, there exists a topology $\mathcal{W}_1$ on $\Ec$, witness of the almost zero-dimensionality of $\Ec$, such that $\Ec$ has a neighborhood basis $\beta_0$ of subsets that are compact in $\mathcal{W}_1$. Let $\mathcal{W}$ be the product topology of $\Q\times\langle\Ec,\mathcal{W}_1\rangle$. Let $\beta$ be the collection of all sets of the form $V\times B$, where $V$ is non-empty and clopen in $\Q$,  and $B\in\beta_0$. Choose a sequence $\{F_n\colon n\in\omega\}$ of compact subsets of $\Q$ such that (i) $F_n\subset F_{n+1}$ for every $n\in\omega$, (ii) $F_{n+1}\setminus F_n$ is countable discrete, and dense in $F_{n+1}$ for every $n\in\omega$, and (iii) $\Q=\bigcup\{F_n\colon n\in\omega\}$.

Let $E_n=F_n\times\Ec$ for every $n\in\omega$. We claim that the topology $\mathcal{W}$, the collection $\{E_n\colon n\in\omega\}$ and $\beta$ satifsy the conditions in Definition \ref{defi:class-sE} for $\Q\times\Ec$.

First, notice that $\mathcal{W}$ witnesses that $\Q\times\Ec$ is almost zero-dimensional. Conditions (a), (b) and (c) follow directly from our choices.

Next, we prove (d). Let $\langle x,y\rangle\in\Q\times\Ec$ and let $m=\min\{k\in\omega\colon x\in F_k\}$. Since $\Ec$ is cohesive, there exists an open set $U$ of $\Ec$ such that $x\in U$ and $U$ contains no non-empty clopen subsets. Let $V$ be open in $\mathbb{Q}$ such that $x\in V$ and $V\cap F_{k}=\emptyset$ if $k<m$. Define $W=V\times U$. Let $n\in\omega$, we argue that $W\cap E_n$ contains no non-empty clopen sets. This is clear if $n<m$ so consider the case when $n\geq m$.  Assume that $O\subset W\cap E_n$ is clopen and non-empty, and consider $\langle a,b\rangle\in O$. Then $(\{a\}\times\Ec)\cap O$ is a non-empty clopen subset of $\{a\}\times\Ec$ such that $(\{a\}\times\Ec)\cap O\subset\{a\}\times U$. This is a contradiction to our choice of $U$. We conclude that (d) holds.

Finally, let us prove (e). Let $V\times B\in\beta$ and $n\in\omega$. Then $(V\times B)\cap E_n=(V\cap F_n)\times B$, which is compact. Also, it is clear that $\beta$ is a basis for the topology of $\Q\times\Ec$. This completes the proof of this result.
\end{proof}

\begin{propo}\label{propo:class-sE-inside-sL}
 If $E\in\sigma\mathcal{E}$ then there exists $\langle C,X,\varphi\rangle\in\sigma\mathcal{L}$ that generates $E$.
\end{propo}
\begin{proof}
From Definition \ref{defi:class-sE}, let us consider for $E$: the witness topology $\mathcal{W}$, the basis $\beta$ of neighborhoods, and the collection $\{E_n\colon n\in\omega\}$. 

We may assume that $\beta$ is countable. For every $B\in\beta$, let $\mathcal{B}_B$ be a countable collection of clopen subsets of $\langle E,\mathcal{W}\rangle$ such that $B=\bigcap\mathcal{B}_B$. Then by a standard Stone space argument, there exists a compact, zero-dimensional and metric space $C$ containing $\langle E,\mathcal{W}\rangle$ as a dense subspace and such that $\cl[C]{O}$ is clopen in $C$ for every $O\in\bigcup\{\mathcal{B}_B\colon B\in\beta\}$. For every $n\in\omega$ let $X_n=\cl[C]{E_n}$; notice that $X_n\cap E=E_n$ since $E_n$ is closed in $\mathcal{W}$. Define $X=\bigcup\{X_n\colon n\in\omega\}$.

We claim that $X$ is witness to the almost zero-dimensionality of $E$; we will prove that $B$ is closed in $X$ for every $B\in\beta$. It is enough to prove that if $m\in\omega$ and $B\in\beta$ are fixed then
$$
\left(\bigcap\{\cl[C]{O}:O\in\mathcal{B}_B\}\right)\cap X_m=B\cap X_m.\ (\ast)
$$
The right side of $(\ast)$ is contained in the left side by the definition of $\mathcal{B}_B$. So take $z\in C$ that is not on the right side of $(\ast)$, we will prove that it is not on the left side.

We may assume that $z\in X_m$. By the choice of $\beta$, we know that $B\cap X_m$ is compact. So there is an open set $U$ of $C$ such that $z\in U$ and $\cl[C]{U}\cap (B\cap X_m)=\emptyset$. Let $F=\cl[C]{U}\cap E_m$. Notice that $F$ is closed in $\langle E_m,\mathcal{W}\restriction E_m\rangle$, and thus, in $\langle E,\mathcal{W}\rangle$. Also, since $U\cap X_n$ is open in $X_n$, $E_n$ is dense in $X_n$ and $z\in U\cap X_n$, then it easily follows that $z\in\cl[C]{F}$. Finally, $F$ is disjoint from $B$ because $F\cap B=(\cl[C]{U}\cap E_m)\cap B=\cl[C]{U}\cap (B\cap E_m)=\cl[C]{U}\cap (B\cap X_m)=\emptyset$. Then $F$ and $B$ are two disjoint closed subsets in $\langle E,\mathcal{W}\rangle$ so there exists $O\in\mathcal{B}_B$ such that $O\cap F=\emptyset$. Since $\cl[C]{O}$ is open in $K$ and disjoint from $F$, it is also disjoint from $\cl[C]{F}$. But $z\in\cl[C]{F}$, so $z\notin\cl[C]{O}$. This shows that $z$ is not on the left side of $(\ast)$.

We have proved that $X$ is witness to the almost zero-dimensionality of $E$. By Lemma 4.11 of \cite{ME} there exists a USC function $\psi_0\colon X\to[0,1)$ such that $\psi_0^\leftarrow(0)=X\setminus E$ and the function $h_0\colon E\to G_0^{\psi_0}$ defined by $h_0(x)=\langle x,\psi_0(x)\rangle$ is a homeomorphism. By condition (d) in Definition \ref{defi:class-sE} we know that $G_0^{\psi_0}$ is $\{G_0^{\psi_0\restriction X_n}\colon n\in\omega\}$-cohesive. Moreover, $\{x\in X_n\colon\psi_0(x)>0\}=E_n$ is dense in $X_n$ for every $n\in\omega$. Lemma 5.9 of \cite{ME} tells us that we can find a USC function $\psi_1\colon X\to[0,1)$ such that $\psi_1\restriction{X_n}$ is a Lelek function for each $n\in\omega$, and the function $h_1\colon G_0^{\psi_0}\to G_0^{\psi_1}$ given by $h_1(\langle x,\psi_0(x)\rangle)=\langle x,\psi_1(x)\rangle$ is a homeomorphism. Now, let $\varphi=\ext{\psi_1}\colon C\to[0,1)$. Then $\langle C,X,\varphi\rangle$ can be easily seen to be an element of $\sigma\mathcal{L}$ and $h_1\circ h_0\colon E\to G_0^{\varphi\restriction X}$ is a homeomorphism.This completes the proof of this result.
\end{proof}

Our main result will be the following.

\begin{theorem}\label{thm-main}
 Let $E$ be a space. Then the following are equivalent:
 \begin{enumerate}[label=(\roman*)]
  \item $E\in\sigma\mathcal{E}$,
  \item there exists $\langle C,X,\varphi\rangle\in\sigma\mathcal{L}$ that generates $E$, and
  \item $E$ is homeomorphic to $\Q\times\Ec$.
 \end{enumerate}
\end{theorem}

The proof of Theorem \ref{thm-main} will be given as follows. First, notice that by Proposition \ref{propo:class-sE-inside-sL}, (i) implies (ii). That (iii) imples (i) is Lemma \ref{lema:QtimesEc-sigmaE}. Also, by Lemma \ref{lema:QtimesEc-sigmaE}, $\sigma\mathcal{E}$ is non-empty so $\sigma\mathcal{L}$ is non-empty as well. Thus, in order to prove that (ii) implies (iii) it is enough to show that any two spaces generated by triples of $\sigma\mathcal{L}$ are homeomorphic. This will be the content of Section \ref{section:main-theorem}.

Given a separable metrizable space $X$, in \cite{vm-cap} $\mathsf{CAP}(X)$ is defined to be the class of separable metrizable spaces $Y=\bigcup\{X_n\colon n\in\omega\}$ such that $X_n$ is closed in $X$, $X_n$ is a nowhere dense subset of $X_{n+1}$ and $X_n\approx X$ for each $n\in\omega$. So $\sigma\mathcal{E}\subset\mathsf{CAP}(\Ec)$ but we do not know whether the other inclusion holds. 

\begin{ques}\label{ques:cap}
 Is $\sigma\mathcal{E}=\mathsf{CAP}(\Ec)$?
\end{ques}

\section{Uniqueness theorem}\label{section:main-theorem}

In this section we give the proof of Theorem \ref{thm-main}. Let $\varphi,\psi\colon X\to[0,\infty)$ be USC functions. In chapter 6 of \cite{ME}, $\varphi$ and $\psi$ are defined to be \emph{$m$-equivalent} if there is a homeomorphism $h\colon X\to Y$ and a continuous function $\alpha\colon X\to(0,\infty)$ such that $\psi\circ h=\alpha\cdot\varphi$. It follows that when $\varphi$ and $\psi$ are $m$-equivalent then $G_0^\varphi$ is homeomorphic to $G_0^\psi$. So, according to the discussion at the end of the previous section, in order to prove Theorem \ref{thm-main}, it is sufficient to prove the following statement.

\begin{propo}\label{propo-equivalence}
 Let $\langle C,X,\varphi\rangle,\langle D,Y,\psi\rangle\in\sigma\mathcal{L}$. Then there exists a homeomorphism $h\colon C\to D$ and a continuous function $\alpha\colon C\to (0,\infty)$ such that $f[X]=Y$ and $\psi\circ h=\alpha\cdot \varphi$.
\end{propo}

The rest of this section will consist on a proof of Proposition \ref{propo-equivalence}. The construction of the homeomorphism $h$ will require us to use two different techniques and mix them. First, we need the tools used in \cite{ME} to extend homeomorphisms using Lelek functions. 

\begin{theorem}\label{thm:Lelek-uniqueness}\cite[Theorem 6.2, p. 26]{ME}
If $\varphi\colon C\to[0,\infty)$ and $\psi\colon D\to[0,\infty)$ are Lelek functions with $C$ and $D$ compact, and $t>\lvert\log{\left(M(\psi)/M(\varphi)\right)}\rvert$, then there exists a homeomorphism $h\colon C\to D$ and a continuous function $\alpha\colon C\to(0,\infty)$ such that $\psi\circ h=\alpha\cdot\varphi$ and $M(\log\circ\alpha)<t$.
\end{theorem}

\begin{theorem}\label{thm:homeo-extension-Lelek}\cite[Theorem 6.4, p. 28]{ME}
Let $\varphi\colon C\to[0,\infty)$ and $\psi\colon D\to[0,\infty)$ be Lelek functions with $C$ and $D$ compact. Let $A\subset C$ and $B\subset D$ be closed such that $G_0^{\varphi\restriction A}$ and $G_0^{\psi\restriction B}$ are nowhere dense in $G_0^\varphi$ and $G_0^\psi$, respectively. Let $h\colon A\to B$ be a homeomorphism and $\alpha\colon A\to(0,\infty)$ a continuous function such that $\psi\circ h=\alpha\cdot{(\varphi\restriction A)}$. If $t\in\mathbb{R}$ is such that $t>\lvert\log{\left(M(\psi)/M(\varphi)\right)}\rvert$ and $M(\log\circ\alpha)<t$ then there is a homeomorphism $H\colon C\to D$ and a continuous function $\beta\colon C\to (0,\infty)$ such that $H\restriction A=h$, $\beta\restriction C=\alpha$, $\psi\circ H=\beta\cdot\varphi$ and $M(\log\circ\beta)<t$.
\end{theorem}

Theorem \ref{thm:Lelek-uniqueness} is called the Uniqueness Theorem for Lelek functions; Theorem \ref{thm:homeo-extension-Lelek} is the Homeomorphism Extension Theorem for Lelek funcions.

The second tool we will need is that of \emph{Knaster-Reichbach covers}. KR-covers were used by Knaster and Reichbach~\cite{KnasterReichbach} to prove homeomorphism extension results in the class of all zero-dimensional spaces. The term KR-cover was first used by van Engelen~\cite{fons86} who proved their existence in a general setting. However, in this paper we will not need the existence of KR-covers in general. We will only need the following straightforward result which is a specific case of KR-covers.


\begin{lemma}\label{lema:KR-restricted}
Fix a metric on $2^\omega$. Let $F\subset 2^\omega$ be closed and assume that $\mathcal{U}=\{U_n:n\in\omega\}$ is a partition of $2^\omega\setminus F$ into clopen sets such that for every $\epsilon>0$ the set $\{n\in\omega:\mathrm{diam}{(U_n)}\geq\epsilon\}$ is finite. Assume that $h\colon 2^\omega\to 2^\omega$ has the following properties
\begin{enumerate}
 \item $h$ is a bijection,
 \item $h\restriction F=\mathsf{id}_{F}$,
 \item for each $n\in\omega$, $h[U_n]=U_n$, and
 \item for each $n\in\omega$, $h\restriction U_n\colon U_n\to U_n$ is a homeomorphism.
\end{enumerate}
Then $h$ is a homeomorphism.
\end{lemma}

We then remark that our proof will be an amalgamation of the Dijkstra-van Mill proof of Theorem 7.5 from \cite{ME} and the van Engelen proof of Theorem 3.2.6 from \cite{fons86}. The functions $h$ and $\alpha$ in the statement of Proposition \ref{propo-equivalence} will be uniform limits of functions. The following discussion can be found in \cite{vm-inf_dim_funct_spaces}. 

Let $X$ and $Y$ be compact metrizable spaces and let $\rho$ be a metric on $Y$. In the set $C(X,Y)=\{f\in Y^X\colon f\,\textrm{is continuous}\}$ we define the uniform metric $\rho$ by $\rho(f,g)=\sup\{\rho(f(x),g(x))\colon x\in X\}$, when $f,g\in C(X,Y)$. It is known that this metric is complete so we may construct complicated continuous functions using Cauchy sequences of simpler continuous functions.

For a compact space $X$, $\mathcal{H}(X)$ denotes the subset of $C(X,X)$ consisting of homeomorphisms. However, even though Cauchy sequences of homeomorphisms will converge to continuous functions, they will not necessarily converge to a homeomorphism. In order to achieve this, we will use the \emph{Inductive Convergence Criterion}. We present the statement of this criterion as it appears in \cite{fons86}.

\begin{theorem}\label{inductive-convergence}\cite[Lemma 3.2.5]{fons86}
 Let $X$ be a zero-dimensional compact metric space with metric $\rho$ and for each $n\in\omega$, let $h_n\colon X\to X$ be a homeomorphism. If for every $n\in\omega$ we have that $\rho(h_{n+1}, h_n)<\epsilon_n $, where
$$\epsilon_n=\min \{2^{-n},3^{-n} \cdot \min\{\min\{\rho(h_i(x),h_i(y))\colon x,y\in X,\,\rho(x,y)\geq 1/n\}: i\leq n \}\},$$
 then the uniform limit $h=\lim_{n\to \infty} h_n$ is a homeomorphism.
\end{theorem}

The exact values of the numbers $\epsilon_n$ in the statement of Theorem \ref{inductive-convergence} are not important. What we will use is that $\epsilon_n$ is a positive number than can be calculated once the first $n+1$ homeomorphisms $h_0,\ldots,h_n$ have been defined.

Before we continue with the proof of Proposition \ref{propo-equivalence}, we stop to give two final ingredients in the proof. 

\begin{lemma}\label{lema:USC-replace-Lelek}
 If $\langle C,X,\psi\rangle\in\sigma\mathcal{L}$ then there exists a Lelek function $\varphi\colon C\to[0,1]$ such that $\langle C,X,\varphi\rangle\in\sigma\mathcal{L}$, $\varphi\restriction X=\psi\restriction X$ and the graph of $\varphi\restriction X$ is dense in the graph of $\varphi$.
\end{lemma}
\begin{proof}
 Let $d_0$ be a metric for $C$ and consider the metric $d(\langle x,y\rangle,\langle z,w\rangle)=d_0(x,z)+\lvert y-w\rvert$ defined on $C\times[0,1]$. Define $\varphi=\ext{\psi\restriction X}$. 
 
 We show that $\varphi$ is a Lelek function. Let $p\in C$ with $\varphi(p)>0$, $t\in(0,\varphi(p))$ and $\epsilon>0$, we want to find $q\in G_0^\varphi$ such that $d(q,\langle p,t\rangle)<\epsilon$. By Lemma \ref{lema:ext-USC} we know that the graph of $\psi\restriction X$ is dense in the graph of $\varphi$ so there exists $k\in\omega$ and $x\in X_k$ such that $d(\langle x,\psi(x)\rangle,\langle p,\varphi(p)\rangle)<\epsilon/2$. We may also assume that $\psi(x)>t$. Since $\psi\restriction X_k$ is a Lelek function, there is $z\in X_k$ such that $d(\langle z,\psi(z)\rangle,\langle x,t\rangle)<\epsilon/2$. So let $q=\langle z,\psi(z)\rangle$. We know that $\psi(z)=\varphi(z)$ so $q\in G_0^\varphi$. Then
 $$
 \begin{array}{lcl}
  d(q,\langle p,t \rangle) & = & d_0(z,p)+\lvert \psi(z)-t\rvert \\  & \leq & d_0(z,x)+d_0(x,p)+\lvert \psi(z)-t\rvert\\
  & = & d(\langle z,\psi(z)\rangle,\langle x,t\rangle)+d_0(x,p)\\
  & \leq & d(\langle z,\psi(z)\rangle,\langle x,t\rangle)+ d(\langle x,\psi(x)\rangle,\langle p,\varphi(p)\rangle)\\
  & < & \epsilon/2+\epsilon/2\\
  & = & \epsilon.
 \end{array} 
 $$
 This shows that $\varphi$ is a Lelek function. The remaining condition holds directly from Lemma \ref{lema:ext-USC}.
\end{proof}

The constant function with value $1$ will be denoted by $\mathsf{1}$.

\begin{lemma}\label{lema:auxiliar-continuous}
 Let $F\subset 2^\omega$ be closed and let $\{V_n\colon n\in \omega\}$ be a partition of $2^\omega\setminus F$ into clopen non-empty subsets. Assume that $\alpha\colon 2^\omega\to (0,\infty)$ has the following properties
\begin{enumerate}
\item $\alpha\restriction F=\mathsf{1}\restriction F$,
\item $lim_{n\to \infty}M({\log\circ(\alpha}\restriction V_n))=0$, and
\item $\alpha\restriction V_n$ is continuous for each $n\in \omega$.
\end{enumerate}
Then $\alpha$ is continuous.
\end{lemma}
\begin{proof}
It is enough to prove that if $\langle x_i\colon i\in\omega\rangle$ is a sequence contained in $2^\omega\setminus F$ such that $x=\lim_{i\to\infty}{x_i}\in F$, then $\lim_{i\to\infty}{\alpha(x_i)}=1$. 

Let $\epsilon>0$. By the continuity of the exponential function there is $\delta>0$ be such that if $t\in(-\delta,\delta)$ then $e^t\in(1-\epsilon,1+\epsilon)$. By condition (2) there exists $N\in\omega$ such that if $n\geq N$, then $\lvert M(\log\circ(\alpha\restriction V_n))\rvert<\delta$. On the other hand, there exists $k\in\omega$ such that if $i>k$ then $x_i\in\bigcup\{V_n\colon n\geq N\}$. If $i\geq k$ we obtain that $\lvert \log{(\alpha(x_i))}\rvert<\delta$ so $\log(\alpha(x_i))\in(-\delta,\delta)$. Thus, $\alpha(x_i)\in(1-\epsilon,1+\epsilon)$ so $\lvert \alpha(x_i)-1\rvert<\epsilon$.
\end{proof}

We now prove our main result. In our proof we will use the tree $\omega^{<\omega}$ of finite sequences of natural numbers. This includes the concatenation $s^\frown i$ where $s\in\omega^{<\omega}$ and $i\in\omega$, that is, the unique sequence with $\mathrm{dom}(s^\frown i)=\mathrm{dom}(s)+1$, $s\subset s^\frown i$ and $(s^\frown i)(\mathrm{dom}(s))=i$.

\begin{proof}[Proof of Proposition \ref{propo-equivalence}]

Without loss of generality we assume that $C=D=2^\omega$, and we fix some metric $\rho$ on $2^\omega$. By an application of Lemma \ref{lema:USC-replace-Lelek} we can assume that $\varphi$ and $\psi$ are Lelek functions, that the graph of $\varphi\restriction X$ is dense in the graph of $\varphi$, and that the graph of $\psi\restriction Y$ is dense in the graph of $\psi$. After this, apply Theorem \ref{thm:Lelek-uniqueness}, so we may assume that $\varphi=\psi$. Then $\langle 2^\omega,X,\varphi\rangle,\langle 2^\omega,Y,\varphi\rangle\in\sigma\mathcal{L}$ so there are collections $\{X_n\colon n\in\omega\}$ and $\{Y_n\colon n\in\omega\}$ that satisfy the conditions in Definition \ref{defi:class-L}. Notice that since the graphs of $\varphi\restriction X$ and $\varphi\restriction Y$ are dense in the graph of $\varphi$ it is easy to see that
\begin{quote}
 $(\ast)$ if $U\subset X$ is open then $$M(\varphi\restriction U)=\sup\{M(\varphi\restriction U\cap X_i)\colon i\in\omega\}=\sup\{M(\varphi\restriction U\cap Y_i)\colon i\in\omega\}.$$
\end{quote}

Given $s\in\omega^{<\omega}$, we construct clopen sets $U_s$ and $V_s$ of $2^\omega$, closed nowhere dense sets $D_s$ and $E_s$ of $X$ and $Y$, respectively, and for every $m\in\omega$ a continuous function $\beta_m\colon 2^\omega\to(0,1)$ and a homeomorphism $h_m\colon 2^\omega\to 2^\omega$. We abreviate the composition $h_n\circ\ldots\circ h_0=f_n$ for all $n\in\omega$. We will use the Inductive Convergence Criterion (Theorem \ref{inductive-convergence}) to make the homeomorphisms converge, so at step $n$ we may calculate the corresponding $\epsilon_n>0$. Our construction will have the following properties.

\begin{enumerate}[label=(\alph*)]
 \item $U_\emptyset=V_\emptyset=2^\omega$.
 \item For each $s\in\omega^{<\omega}$, $D_s\subset U_s$ and $E_s\subset V_s$.
 \item For every $n\in\omega$ and $s\in\omega^n$, $\{U_{s^\frown i}\colon i\in\omega\}$ is a partition of $U_s\setminus D_s$ and $\{V_{s^\frown i}\colon i\in\omega\}$ is a partition of $V_s\setminus E_s$.
 \item For every $n\in\omega$, $X_n\subset\bigcup\{D_s\colon s\in\omega^{\leq n}\}$ and $Y_n\subset\bigcup\{E_s\colon s\in\omega^{\leq n}\}$.
 \item For every $n\in\omega$ and $s\in\omega^{n+1}$, $\mathrm{diam}(U_s)\leq 2^{-n}$ and $\mathrm{diam}(V_s)\leq\min\{2^{-n},\epsilon_n\}$.
 \item For every $n\in\omega$ and $s\in\omega^n$, $f_n[D_s]=E_s$.
 \item For every $n\in\omega$ and $s\in\omega^n$, $h_{n+1}\restriction{E_s}=\mathsf{id}_{E_s}$.
 \item For every $n\in\omega$ and $s\in\omega^{n+1}$, $f_{n}[U_s]=V_s$.
 \item For every $n,k\in\omega$, $\{s\in\omega^n\colon \textrm{diam}(U_s)\geq 2^{-k}\}$ is finite.
 \item For every $n\in\omega$ and $x\in 2^\omega$, $\lvert \log(\beta_{n+1}(x)/\beta_n(x))\rvert<2^{-n}$.
 \item For every $n\in\omega$, $\varphi=(\beta_n\cdot\varphi)\circ{f_n^{-1}}$.
\end{enumerate}

Let us assume that we have finished this construction, we claim that $f=\lim_{n\rightarrow\infty}f_n$ exists, is a homeomorphism and $f[X]=Y$. 

First, let $x\in 2^\omega$ and $n\in\omega$. If $x\in\bigcup_{s\in\omega^n}{D_s}$, then $f_n(x)=f_{n+1}(x)$ by conditions (f) and (g). Thus, $\rho(f_n(x),f_{n+1}(x))=0$. Otherwise, by (c) there exists $t\in\omega^{n+1}$ with $x\in U_t$. By (h), $f_n(x)\in V_t$. If $x\in D_t$, by (f) and (b), $f_{n+1}(x)\in E_t\subset V_t$. Otherwise, by (c), there is $i\in\omega$ with $x\in U_{t^\frown i}$ so by (h), $f_{n+1}(x)\in V_{t^\frown i}\subset V_t$. In any case, we obtain that $f_{n+1}(x)\in V_t$. So $\rho(f_n(x),f_{n+1}(x))<\epsilon_n$ by the second part of (e). Thus, $\rho(f_n,f_{n+1})<\epsilon_n$ and we can apply the Inductive Convergence Criterion to conclude that $f$ is well-defined and in fact, a homeomorphism.

Next, let $x\in X$ so $x\in X_m$ for some $m\in\omega$. Thus, by (d) there exists $s\in\omega^{\leq m}$ such that $x\in D_s$. Then $f_{\mathrm{dom}(s)}(x)\in E_s\subset Y$ by (f). By (g) it inductively follows that $f_n(x)=f_{\mathrm{dom}(s)}(x)$ for every $n\geq\mathrm{dom}(s)$. This implies that $f(x)\in Y$. A completely analogous argument shows that if $y\in Y$ then there is $x\in X$ such that $f(x)=y$. This shows that $f[X]=Y$.

By (j) we know that $\{\beta_n\colon n\in\omega\}$ is a Cauchy sequence with the uniform metric so $\beta=\lim_{n\rightarrow\infty}\beta_n$ exists and is a continuous function. Using the first part of (e) it is possible to prove that $\{f_n^{-1}\colon n\in\omega\}$ is also a Cauchy sequence and converges to $f^{-1}$; this proof is completely analogous to the proof that $f=\lim_{n\rightarrow\infty}f_n$ so we omit it. Then, by uniform continuity we infer that $\lim_{n\rightarrow\infty}\beta_n\circ f_n^{-1}=\beta\circ f^{-1}$. So using that $\varphi$ is USC and (k) we obtain the following
$$
\begin{array}{lcl}
 \beta(x)\cdot\varphi(x) & = & \lim_{n\rightarrow\infty}\beta_n(x)\cdot\varphi(x)\\ & = & \lim_{n\rightarrow\infty}\varphi(f_n(x))\\ & \leq &  \varphi(f(x))\\
 & = & \lim_{n\rightarrow\infty}\varphi(f_n(f_n^{-1}(f(x))))\\ & = & \lim_{n\rightarrow\infty}\beta_n(f_n^{-1}(f(x)))\cdot \varphi(f_n^{-1}(f(x)))\\ & \leq & \beta(x)\cdot\varphi(x)
\end{array}
$$

Thus, $\varphi\circ f=\beta\cdot\varphi$. This argument is completely analogous to the one in \cite[Theorem 7.5]{ME}.

Now we carry out the construction. Let $\gamma\colon \omega^{<\omega}\setminus\{\emptyset\}\to\omega$ be any function such that $\gamma\restriction\omega^{m+1}$ is injective for all $m\in\omega$.

\underline{Step $0$}. Let $U_\emptyset=V_\emptyset=2^\omega$, as in condition (a). From $(\ast)$ we infer that there exists $k_\emptyset\in\omega$ such that
$$
\begin{array}{ll}
 \log{(M(\varphi))}-\log{(M(\varphi\restriction X_{k_\emptyset}))}<1, & \textrm{ and} \\[0.5em]
 \log{(M(\varphi))}-\log{(M(\varphi\restriction Y_{k_\emptyset}))}<1. &
\end{array}
$$
Define $D_\emptyset=X_{k_\emptyset}$ and $E_\emptyset=Y_{k_\emptyset}$. Then $\varphi\restriction D_\emptyset$ and $\varphi\restriction E_\emptyset$ are Lelek functions, and $\lvert\log(M(\varphi\restriction E_\emptyset)/M(\varphi\restriction D_\emptyset))\rvert<1$ so we may apply Theorem \ref{thm:Lelek-uniqueness} to obtain a homeomorphism $\widehat{h}_\emptyset\colon D_\emptyset\to E_\emptyset$ and a continuous function $\alpha_\emptyset\colon D_\emptyset\to(0,\infty)$ such that $\varphi\circ\widehat{h}_\emptyset=(\varphi\restriction D_\emptyset)\cdot\alpha_\emptyset$ and $M(\log\circ\alpha_\emptyset)<1$. After this, apply Theorem \ref{thm:homeo-extension-Lelek} to find a homeomorphism $h_0\colon 2^\omega\to 2^\omega$ and a continuous function $\beta_0\colon 2^\omega\to(0,\infty)$ such that $h_0\restriction D_\emptyset=\widehat{h}_\emptyset$, $\beta_0\restriction D_\emptyset=\alpha_\emptyset$, $\varphi\circ h_0=\varphi\cdot\beta_0$ and $M(\log\circ\alpha_0)<1$.

Notice that since $h_0=f_0$ this implies (k) for $n=0$. Let $\{V_n\colon n\in\omega\}$ be a partition of $E_\emptyset$ into clopen sets with their diameters converging to $0$. We may assume that $\mathrm{diam}(V_n)<\min\{\epsilon_0,1\}$ for every $n\in\omega$. We define $U_n=h_0^\leftarrow[V_n]$ for each $n\in\omega$. Without loss of generality we may assume that for all $n\in\omega$, $\mathrm{diam}(U_n)<1$. With this we have finished step $0$ in the construction.

\underline{Inductive step}: Assume that we have constructed the sets $D_s,E_s$ for $s\in\omega^{\leq m}$, the sets $U_s,V_s$ for $s\in\omega^{\leq m+1}$ the homeomorphisms $h_i$ for $i\leq m$, and the continuous functions $\beta_i$ for $i\leq m$. Notice that by condition (c) it inductively follows that $\bigcup\{D_s\colon s\in\omega^{\leq m}\}$ and $\bigcup\{E_s\colon s\in\omega^{\leq m}\}$ are closed because their complement is $\bigcup\{U_s:s\in\omega^{m+1}\}$, and $\bigcup\{V_s:s\in\omega^{m+1}\}$, respectively.

Fix $t\in\omega^{m+1}$. First, notice that by $(\ast)$ we have that there exists $k_t\in\omega$ such that
$$
\log{(M(\varphi\restriction V_t))}-\log{(M(\varphi\restriction V_t\cap Y_{k_t}))}<2^{-(m+\gamma(t))}.
$$
Notice that $\varphi\restriction V_t\cap Y_{k_t}$ is a Lelek function.

Recall that (k) says that $\varphi=(\beta_m\cdot \varphi)\circ f_n^{-1}$. In particular this implies that $\varphi\restriction V_t=(\beta_m\cdot \varphi)\restriction U_t\circ f_n^{-1}\restriction V_t$; from this we infer the following. First, using $(\ast)$ we may assume that $k_t\in\omega$ is such that
$$
 \log{(M(\varphi\restriction V_t))}-\log{(M(\varphi\restriction V_t\cap f_m[X_{k_t}]))}<2^{-(m+\gamma(t))}.
$$
Also, $\varphi\restriction V_t\cap f_m[X_{k_t}]$ is a Lelek function. 

So define $D_t=V_t\cap f_m[X_{k_t}]$ and $E_t=V_t\cap Y_{k_t}$. Then $\varphi\restriction D_t$ and $\varphi\restriction E_t$ are Lelek functions, and $\lvert\log(M(\varphi\restriction E_t)/M(\varphi\restriction D_t))\rvert<2^{-(m+\gamma(t))}$ so we may apply Theorem \ref{thm:Lelek-uniqueness} to obtain a homeomorphism $\widehat{h}_t\colon D_t\to E_t$ and a continuous function $\widehat{\alpha}_t\colon D_t\to (0,\infty)$ such that $\varphi\circ \widehat{h}_t=\varphi\cdot\widehat{\alpha}_t$ and $M(\log\circ\widehat{\alpha}_t)<2^{-(m+\gamma(t))}$. Then apply Theorem \ref{thm:homeo-extension-Lelek} to find a homeomorphism $h_t\colon  V_t\to V_t$ and a continuous function $\alpha_t\colon V_t\to(0,\infty)$ such that $h_t\restriction D_t=\widehat{h}_t$, $\alpha_t\restriction D_t=\widehat{t}_t$, $\varphi\circ h_t=\varphi\restriction V_t\cdot\alpha_t$ and $M(\log\circ\alpha_t)<2^{-(m+\gamma(t))}$.

Let $E_m=\bigcup\{E_s\colon s\in\omega^{\leq m}\}$. Then define
$$
h_{m+1}=\mathsf{id}_{E_m}\cup\bigcup\{h_s\colon s\in\omega^{m+1}\},
$$
by Lemma \ref{lema:KR-restricted} it follows that $h_{m+1}$ is a homeomorphism. Also, define 
$$
\alpha_{m+1}=\mathsf{1}\restriction{E_m}\cup\bigcup\{\alpha_s\colon s\in\omega^{m+1}\},
$$
and $\beta_{m+1}(x)=\alpha_{m+1}(f_m(x))\cdot\beta_m(x)$ for all $x\in 2^\omega$.  By Lemma \ref{lema:auxiliar-continuous}, $\alpha_{m+1}$ is continuous so $\beta_{m+1}$ is continuous.

Now, fix $t\in\omega^{m+1}$ again. Write $V_t\setminus E_t$ as a union of a countable, pairwise disjoint collection of clopen sets, all diameters of which are smaller than $\min\{\epsilon_m,2^{-m}\}$ and converge to $0$. Let $\{V_{t^\frown i}\colon i\in\omega\}$ be such partition and for each $i\in\omega$, let $U_{t^\frown i}=f_{m+1}^{\leftarrow}[V_{t^\frown i}]$. Without loss of generality we may assume that for $i\in\omega$, $\mathrm{diam}(U_{t^\frown i})<2^{-m}$.

We leave the verification that all conditions (a) to (k) hold in this step of the induction to the reader. This concludes the inductive step, and the proof of this result.
\end{proof}

\section{The hyperspace of finite sets of \texorpdfstring{$\Ec$}{Ec}}\label{section:hyperspace}

For a space $X$, $\mathcal{K}(X)$ denotes the hyperspace of non-empty compact subsets of $X$ with the
Vietoris topology. For any $n\in  \N$, $\mathcal{F}_n(X)$ 
is the subspace of $\mathcal{K}(X)$ consisting
of all non-empty subsets that have cardinality less than or equal to $n$, and $\mathcal{F}(X)$ is the subspace of $\mathcal{K}(X)$ of finite non-empty subsets of $X$. 

Given $n\in \N$ and subsets $U_0,\ldots, U_n$ of a topological space $X$, $\vietoris{ U_{0},\ldots ,U_{n}}$ denotes the collection $\left\lbrace   F \in \mathcal{K}(X)\colon F\subset \bigcup_{k=0}^n U_k,\, F\cap U_{k}\neq \emptyset \textit{ for } k \leq n \right\rbrace $. Recall that the Vietoris topology on $\mathcal{K}(X)$ has as its canonical basis all the sets of the form $\vietoris{U_{0},\ldots ,U_{n}}$, where $U_k$ is a non-empty open subset of $X$ for each $k\leq n$.

For each $n\in\N$, let $\pi_n\colon X^n\to \mathcal{F}_n(X)$ be defined by $\pi_n(x_0,\ldots,x_{n-1})=\{x_0,\ldots, x_{n-1}\}$. It is know that this function is continuous, finite-to-one and in fact it is a quotient \cite[2.4.3]{Sub}. 

\begin{lemma}\label{Coh}\cite{zaragoza-2}
Let $X$ be a space that is $\{A_s: s\in S\}$-cohesive, witnessed by a basis $\mathcal{B}$ of open sets. Consider the following collection of subsets of $\mathcal{F}(X)$:
$$
\mathcal{A}=\left\{\pi_n[A_{s_1}\times\cdots\times A_{s_n}]\colon n\in\mathbb{N},\, \forall i\in\{1,\ldots,n\}\, (s_i\in S)\}\right.
$$
Then $\mathcal{F}(X)$ is $\mathcal{A}$-cohesive, and the open sets that witness this may be taken from the collection $\mathcal{C}=\left\{\vietoris{U_1,\ldots, U_n}\colon\forall i\in\{1,\ldots,n\}\, (U_i\in\mathcal{B})\right\}$.
\end{lemma}

Before starting the proof, we remind the reader that if $X$ is separable and metrizable then $\mathcal{K}(X)$ is also separable and metrizable (see \cite[4.5.2]{Sub} and \cite[4.9.13]{Sub}). Thus, with the Vietoris topology we are not leaving our self-imposed universe of discourse.

\begin{propo}
 $\mathcal{F}(\Ec)\in\sigma\mathcal{E}$ 
\end{propo}
\begin{proof}
According to (2) in \cite[Theorem 3.1]{DvM} there is a witness topology $\mathcal{W}_0$ for $\Ec$ and a basis $\beta_0$ for $\Ec$ of sets that are compact in $\mathcal{W}_0$. Let $\mathcal{W}_1$ the Vietoris topology in $\mathcal{K}(\Ec,\mathcal{W}_0)$ and define $\mathcal{W}=\mathcal{W}_1\restriction\mathcal{F}(\Ec)$. Let $\beta$ be the collection of all sets of the form $\vietoris{U_0,\ldots,U_n}\cap\mathcal{F}(\Ec)$ where $n\in\omega$ and $U_j\in\beta_0$ for each $j\leq n$. Also, for every $n\in\omega$ let $E_n=\mathcal{F}_{n+1}(\Ec)$. We will now check that these choices satisfy the conditions in Definition \ref{defi:class-sE}.

By \cite[4.13.1]{Sub} we know that $\mathcal{W}_1$ is zero-dimensional so $\mathcal{W}$ is also zero-dimensional. In \cite[Proposition 2.2]{zaragoza-1} it was proved that $\mathcal{W}$ witnesses that $\mathcal{F}(\Ec)$ is almost zero-dimensional. Condition (a) clearly holds.

For (b), fix $n\in\omega$. Since $\Ec$ is crowded and $\mathcal{F}_{n+1}(\Ec)$ is a continuous image of $\Ec^{n+1}$ (under the function $\pi_{n+1}$ defined above), then $\mathcal{F}_{n+1}(\Ec)$ is crowded. Recall that $\mathcal{F}_n(X)$ is always closed in $\mathcal{K}(X)$ for any topological space $X$ and all $n\in\mathbb{N}$ (\cite[2.4.2]{Sub}). Thus, we only need to show that $\mathcal{F}_{n+2}(\Ec)\setminus\mathcal{F}_{n+1}(\Ec)$ is dense in $\mathcal{F}_{n+2}(\Ec)$; this is well-known but for the reader's convenience we give a short proof. Since $\Ec$ has no isolated points then the set $D$ of all $x\in\Ec^{n+2}$ such that if $i,j\leq n+2$ and $i\neq j$, then $x(i)\neq x(j)$ is easily seen to be dense in $\Ec^{n+2}$. Then $\pi_{n+2}[D]=\mathcal{F}_{n+2}(\Ec)\setminus\mathcal{F}_{n+1}(\Ec)$  is dense in $\mathcal{F}_{n+2}(\Ec)$. This proves (b).

Also, $\mathcal{F}_{n+1}(\Ec)$ is $\mathcal{W}$-closed in $\mathcal{F}(\Ec)$ for all $n\in\omega$, which implies (c). Let $S=\{0\}$ and $A_0=\Ec$. The collection $\mathcal{A}$ from Lemma \ref{Coh} is equal to $\{\mathcal{F}_{n+1}(\Ec)\colon n\in\omega\}$. Thus, by Lemma \ref{Coh} we obtain (d). Finally, it was proved \cite[Proposition 3.4]{zaragoza-1} that if $\mathcal{U}\in\beta$ and $n\in\omega$, then $\mathcal{U}\cap\mathcal{F}_{n+1}(\Ec)$ is compact in $\mathcal{W}\restriction\mathcal{F}_{n+1}(\Ec)$, which implies (e).
\end{proof}

\begin{cor}\label{hyperspace}
 $\mathcal{F}(\Ec)\approx\Q\times\Ec$.
\end{cor}

Here it is natural to ask about $\mathcal{F}(\Q\times\Ec)$, we will prove that this space is homemorphic to $\Q\times\Ec$ as well.

\begin{propo}\label{propo:Fn-E-in-sE}
 Let $E\in\sigma\mathcal{E}$. If $n\in\N$ then $\mathcal{F}_n(E)\in\sigma\mathcal{E}$.
\end{propo}
\begin{proof}
Let $\mathcal{W}$, $\{E_n\colon n\in\omega\}$ and $\beta$ be witnesses of  $E\in\sigma\mathcal{E}$. By \cite[Proposition 2.2]{zaragoza-1}, the Vietoris topology $\mathcal{W}_0$ of $\mathcal{F}_n(E,\mathcal{W})$ witnesses the almost zero-dimensionality of $\mathcal{F}_n(E)$. For each $m\in\omega$, let $Z_m=\pi_m[E_m^n]$. We define $\beta_0$ to be the collection of the sets of the form $\vietoris{U_0,\ldots,U_k}$ where $k<\omega$ and $U_i\in\beta$ for every $i\leq k$. We claim that $\mathcal{W}_0$, $\{Z_m:m\in\omega\}$ and $\beta_0$ witness that $\mathcal{F}_n(E)\in\sigma\mathcal{E}$.

Conditions (a), (b) and (c) are easily seen to follow. By Lemma \ref{Coh}, we infer that $\mathcal{F}_n(E)$ is $\{\mathcal{F}_n(E_m)\colon m\in\omega\}$-cohesive, which is (d). Now, let $U=\vietoris{U_0,\ldots,U_k}\in\beta_0$ and $m\in\omega$. Notice that $U\cap Z_m\subset\vietoris{U_0\cap E_m,\ldots, U_k\cap E_m}$. Now, by the choice of $\beta$ we know that $U_i\cap E_m$ is compact in $\mathcal{W}$ for every $i\leq k$. Thus, the set $\vietoris{U_0\cap E_m,\ldots, U_k\cap E_m}$ is compact in $\mathcal{W}_0$. Since $U\cap Z_m$ is closed in $\mathcal{W}_0$, it is also compact. This proves (e) and completes the proof.
\end{proof}

\begin{propo}\label{propo:F-E-in-sE}
 If $E\in\sigma\mathcal{E}$, then $\mathcal{F}(E)\in\sigma\mathcal{E}$.
\end{propo}
\begin{proof}
 Let $\mathcal{W}$, $\{E_n\colon n\in\omega\}$ and $\beta$ be witnesses of  $E\in\sigma\mathcal{E}$. Let $\mathcal{W}_0$ be the Vietoris topology of $\mathcal{F}(E,\mathcal{W})$. For each $m\in\omega$, let $Z_m=\pi_n[E_m^m]$. We define $\beta_0$ to be the collection of the sets of the form $\vietoris{U_0,\ldots,U_k}$ where $k<\omega$ and $U_i\in\beta$ for every $i\leq k$. The proof that $\mathcal{W}_0$, $\{Z_m:m\in\omega\}$ and $\beta_0$ witness that $\mathcal{F}(E)\in\sigma\mathcal{E}$ is completely analogous to the proof of Proposition \ref{propo:Fn-E-in-sE} and we will leave it to the reader.
\end{proof}

\begin{cor}\label{hyperspace-extra}
 If $n\in\N$, then $\mathcal{F}_n(\Q\times\Ec)\approx\Q\times\Ec$. Also, $\mathcal{F}(\Q\times\Ec)\approx\Q\times\Ec$.
\end{cor}

\section{The \texorpdfstring{$\sigma$}{s}-product of \texorpdfstring{$\Ec$}{Ec}}\label{section:sigma-product}

Given a space $X$, a cardinal $\kappa$ and $e\in X$, the \emph{support} of $x$ with respect to $e$ is the set $\mathrm{supp}_{e}(x)=\{\alpha\in \kappa\colon x(\alpha)\neq e\}$. Then the $\sigma$-product of $\kappa$ copies of $X$ with basic point $e$ is $\sigma(X,e)^\kappa=\{x\in X^\kappa \colon\lvert \mathrm{supp}_{e}(x)\rvert <\omega\}$ as a subspace of $X^\kappa$.
It is known that $\sigma(X,e)^\kappa$ is dense in $X^\kappa$. 

Now, consider $X=\Ec$. Since $\Ec$ is homogeneous, the choice of the point $e$ is irrelevant. Denote $\sigma(\Ec^\omega,e)=\sigma\Ec^\omega$. Since $\sigma{}\Ec^\omega$ is separable and metrizable, it is natural to ask the following.

\begin{ques}\label{ques:sigma}
 Is $\sigma{\Ec}^\omega$ homeomorphic to $\Q\times\Ec$?
\end{ques}

We were unable to answer this question but we make some comments. At first, it seems that it would be possible to prove that $\sigma{\Ec}^\omega\in\sigma\mathcal{E}$ using the following stratification. Given $n\in\omega$, define $\sigma_n\Ec=\{x\in \Ec^\omega\colon\mathrm{supp}_e(x)\subset n\}$. 
It is easy to see that $\sigma_n\Ec$ is closed in $\Ec^\omega$ and homeomorphic to $\Ec^n$ for each $n\in\omega$; so in fact it is a closed copy of $\Ec$ if $n\neq 0$. In fact, using an argument similar to the one in  \cite[Remark 5.2, p. 21]{ME} it is possible to prove the following.

\begin{lemma}\label{lema:sigma-cohesive}
 $\sigma{\Ec}^\omega$ is $\{\sigma_n{\Ec}\colon n\in \N\}$-cohesive.
\end{lemma}
%

Also, a natural witness topology for $\sigma{\Ec^\omega}$ can be obtained by using the restriction of the product topology of the witness topology for $\Ec$. The reader will not find it difficult to prove that properties (a) to (d) of Definition \ref{defi:class-sE} hold but property (e) does not hold. Thus, it is possible that $\sigma\Ec^\omega$ is a different type of space from $\Q\times\Ec$. Notice that a negative answer to Question \ref{ques:sigma} implies a negative answer to Question \ref{ques:cap}.

\section{Factors of \texorpdfstring{$\Q\times\Ec$}{QxEc}}\label{section:factors}

Recall that a space $X$ is a factor of a space $Y$ if there is another space $Z$ such that $X\times Z\approx Y$. In \cite{DvM} the factors of $\Ec$ were characterized and in \cite{ME} the factors of $\E$ were characterized. So we found it natural to try to characterize the factors of $\Q\times\Ec$.

\begin{lemma}\label{lema:copies-inside-QEc}{\phantom{blank space}}
\begin{enumerate}[label=(\alph*)]
 \item $\Q\times\Ec$ does not contain any closed subspace homeomorphic to $\Ec^\omega$.
 \item $\Q\times\Ec$ does not contain any closed subspace homeomorphic to $\E$.
\end{enumerate}
\end{lemma}
\begin{proof}
Assume that $e\colon \Ec^\omega\to\Q\times\Ec$ is a closed embedding. Choose some enumeration $\Q=\{q_n\colon n\in\omega\}$. Notice that $F_n=e^\leftarrow[\{q_n\}\times \Ec]$ is a closed subset of $\Ec^\omega$ for every $n\in\omega$. By the Baire category theorem there exists $m\in\omega$ such that $F_m$ has non-empty interior in $\Ec^\omega$. Recall that every open subset of $\Ec^\omega$ has a closed copy of $\Ec^\omega$ (see the proof of \cite[Corollary 3.2]{d-vm-s}). Thus, this implies that there is a closed copy of $\Ec^\omega$ in $\{q_m\}\times\Ec$. However, $\Ec^\omega$ is cohesive by \cite[Remark 5.2]{ME} and every closed cohesive subset of $\Ec$ is homeomorphic to $\Ec$ by \cite[Theorem 3.5]{DvM}.  This is a contradiction to \cite[Corollary 3.2]{d-vm-s}. Thus, (a) holds.

Now, assume that $e\colon\E\to\Q\times\Ec$ is a closed embedding. Again, let $\Q=\{q_n\colon n\in\omega\}$ be an enumeration and let $F_n=e^\leftarrow[\{q_n\}\times\Ec]$ for every $n\in\omega$. Since $e$ is a closed embedding, for every $n\in\omega$, $F_n$ is homeomorphic to a closed subset of $\Ec$ so it is completely metrizable. This implies that $\E$ is an absolute $G_{\delta\sigma}$, and this contradicts \cite[Remark 5.5]{ME}. This completes the proof of (b).
\end{proof}

\begin{theorem}\label{thm:factors}
 For a non-empty space $E$ the following are equivalent:
 \begin{enumerate}[label=(\arabic*)]
  \item $E\times(\Q\times\Ec)$ is homeomorphic to $\Q\times\Ec$,
  \item $E$ is a $(\Q\times\Ec)$-factor,
  \item there are a topology $\mathcal{W}$ on $E$ witnessing that $E$ is almost zero-dimensional, a collection of $\mathcal{W}$-closed non-empty subsets $\{E_n:n\in\omega\}$ and a basis of neighborhoods $\beta$ such that  \begin{enumerate}[label=(\roman*)]
   \item $E=\bigcup\{E_n\colon n\in\omega\}$,
   \item for every $n\in\omega$, $E_n\subset E_{n+1}$, and
   \item for every $U\in\beta$ and $n\in\omega$, $U\cap E_n$ is compact in $\mathcal{W}$.
  \end{enumerate}
 \end{enumerate}
\end{theorem}
\begin{proof}
 Condition (1) clearly implies (2). 
 
 Next, we prove that (2) implies (3). Since $E$ is a $\Q\times\Ec$-factor, there is a space $Z$ such that $E\times Z\approx\Q\times\Ec$. Let $\mathcal{W}$, $\{X_n\colon n\in\omega\}$ and $\beta$ be witnesses of $E\times Z\in\sigma\mathcal{E}$ as in Definition \ref{defi:class-sE}. Fix $a\in Z$ and let $A=E\times\{a\}$; we may choose $a$ in such a way that $A\cap E_0\neq\emptyset$. We define $E_n=X_n\cap A$ for every $n\in\omega$, $\mathcal{W}_0=\mathcal{W}\restriction A$ and $\beta_0=\{U\cap A\colon U\in\beta\}$. It is not hard to prove that these sets have the corresponding properties (i), (ii) and (iii) replacing $E$ for $A$.
 
 Finally, we prove that (3) implies (1). Let $\mathcal{W}_0$, $\{E_n\colon n\in\omega\}$ and $\beta_0$ as in item (3) for $E$. Let $\mathcal{W}$, $\{X_n\colon n\in\omega\}$ and $\beta$ witnessing that $\Q\times\Ec$, as in Lemma \ref{lema:QtimesEc-sigmaE}. Let $\mathcal{W}_1$ be the product topology of $\langle E, \mathcal{W}_0\rangle\times\langle \Q\times\Ec,\mathcal{W}\rangle$. Notice that $E_n\times X_n$ is $\mathcal{W}_1$-closed for every $n\in\omega$.  Thus, $\mathcal{W}_1$ clearly witnesses that $E\times(\Q\times\Ec)$ is almost zero-dimensional. Finally, let $\beta_1=\{U\times V\colon U\in\beta_0,V\in\beta_1\}$. 
 
 We claim that $\mathcal{W}_1$, $\{E_n\times X_n\colon n\in\omega\}$ and $\beta_1$ witness that $E\times(\Q\times\Ec)\in\sigma\mathcal{E}$. Conditions (a), (b) and (c) are easily checked. By \cite[Remark 5.2]{ME} we obtain that $E\times(\Q\times\Ec)$ is $\{E_n\times X_n\colon n\in\omega\}$-cohesive. Finally, given $U\times V\in\beta_1$ and $n\in\omega$, since $U\cap E_n$ is compact in $\mathcal{W}_0$ and $V\cap X_n$ is compact in $\mathcal{W}$, then $(U\times V)\cap(E_n\times X_n)$ is compact in $\mathcal{W}_1$. This concludes the proof. 
\end{proof}

\begin{ques}\label{ques:factor}
 Can we remove mention of the zero-dimensional witness topology in Theorem \ref{thm:factors} by adding the following statement?\\[0.5em]
  (4) $E$ is a union of a countable collection of $C$-sets, each of which is a $\Ec$-factor.
\end{ques}

David Lipham has informed us that, however, if we change ``$C$-sets'' to ``closed sets'' in (4) of Question \ref{ques:factor}, the resulting statement is not equivalent to $E$ being an $(\Q\times\Ec)$-factor. This is because in \cite{lipham} he gave an example of an $F_\sigma$ subset of $\Ec$ that is not an $\E$-factor.

\begin{cor}
 \begin{enumerate}[label=(\roman*)]
  \item Every $\Ec$-factor is a $(\Q\times\Ec)$-factor.
  \item The space $\Q$ is a $(\Q\times\Ec)$-factor but is not a $\Ec$-factor.
  \item Every $(\Q\times\Ec)$-factor is a $\E$-factor.
  \item The space $\E$ is a $\E$-factor that is not a $(\Q\times\Ec)$-factor.
  \item The space $\Ec^\omega$ is a $\E$-factor that is not a $(\Q\times\Ec)$-factor.
 \end{enumerate}
\end{cor}
\begin{proof}
For (i), let $X$ be a $\Ec$-factor. By \cite[3.2]{DvM}, $X\times\Ec\approx\Ec$. Thus, $X\times(\Q\times\Ec)\approx\Q\times(X\times\Ec)\approx \Q\times\Ec$. For (ii), notice that since $\Q\times\Q\approx\Q$ then $\Q$ is a $(\Q\times\Ec)$-factor but it is not a $\Ec$-factor because it is not Polish. For (iii), let $X$ be a $(\Q\times\Ec)$-factor. By \cite[Proposition 9.1]{ME}, $\Ec\times\Q^\omega\approx\E$. Thus, $X\times\E\approx X\times(\Ec\times\Q^\omega)\approx X\times(\Q\times\Ec)\times\Q^\omega\approx (\Q\times\Ec)\times\Q^\omega\approx\Ec\times\Q^\omega\approx\E$. For (iv), it is clear that $\E$ is a $\E$-factor. 
However, $\E$ is not a $(\Q\times\Ec)$-factor because in that case $\Q\times\Ec$ would have a closed copy of $\E$ and we have proved that this is impossible in Lemma \ref{lema:copies-inside-QEc}. For (v), recall that $\Ec^\omega$ is an $\E$-factor by \cite[Corollary 9.3]{ME} and it cannot be a $(\Q\times\Ec)$-factor, again by Lemma \ref{lema:copies-inside-QEc}.
\end{proof}

\section{Dense embeddings of \texorpdfstring{$\Q\times\Ec$}{QxEc}}\label{section:embeddings}

In this section we consider when $\Q\times\Ec$ can be embedded in almost zero-dimensional spaces as dense subsets. Since every countable dense subset of $\Ec$ is homeomorphic to $\Q$ and $\Ec^2\approx\Ec$, we obtain the following.

\begin{ex}
 There is a dense $F_\sigma$ subset of $\Ec$ that is homeomorphic to $\Q\times\Ec$.
\end{ex}

Moreover, using an analogous argument, $\Ec^\omega$ can be shown to contain dense subsets that are homeomorphic to $\Q\times\Ec^\omega$ so they are non-homeomorphic to $\Q\times\Ec$ by Lemma \ref{lema:copies-inside-QEc}. Thus, we make the following questions.

\begin{ques}\label{ques:dense-in-Ec}
Let $X\subset\Ec$ be dense and a countable union of nowhere dense $C$-sets. If $X$ is cohesive, is it homeomorphic to $\Q\times\Ec$?
\end{ques}

\begin{ques}\label{ques:copy-inside-stable}
 Is there a dense $F_\sigma$ subset of $\Ec^\omega$ that is homeomorphic to $\Q\times\Ec$?
\end{ques}

Notice that Question \ref{ques:dense-in-Ec} is related to Question \ref{ques:cap}. Also, a positive answer to Question \ref{ques:sigma} implies a positive answer to Question \ref{ques:copy-inside-stable}. We recall that it is still unkown whether the hyperspace $\mathcal{K}(\Ec)$ is homeomorphic to $\Ec$ or $\Ec^\omega$ (see Question 5.5 of \cite{zaragoza-1}) but now we know that it has a dense copy of $\Q\times\Ec$ by Corollary \ref{hyperspace}.

\end{document}